\documentclass[10pt]{amsart}
\usepackage{import}
\usepackage{./preamble/Pakete}
\usepackage{./preamble/Befehle}
\begin{document}

\title[$K$-motives and Springer Theory]{K-motives, Springer Theory and the Local Langlands Correspondence}

\author{Jens Niklas Eberhardt}
\email{mail@jenseberhardt.com}  
\begin{abstract}
We construct a geometric realization of categories of representations of affine Hecke algebras and split reductive $p$-adic groups via a $K$-motivic Springer theory.
We suggest a connection to the coherent Springer theory of Ben-Zvi, Chen, Helm, and Nadler through a categorical Chern character and outline results and conjectures on $K$-motives within the Langlands program.

To achieve our results, we introduce a six functor formalism for reduced $K$-motives applicable to linearly reductive stacks and establish formality for categories of Springer $K$-motives. We work within a broader framework of Hecke algebras derived from Springer data. This makes the results applicable, for example, to the ($K$-theoretic) quiver Hecke and Schur algebra. Moreover, we relate our constructions to prior geometric realizations for graded Hecke algebras.

\end{abstract}
\maketitle
\setcounter{tocdepth}{1} 
\tableofcontents

\section{Introduction}
\subsection{Motivation}
Consider a reductive group $G/k$, with a simply connected derived subgroup, defined over an algebraically closed field $k$ of characteristic zero. Denote by $\Hh_{\op{aff},\mathbf{q}}$ the generic affine Hecke algebra associated to $G.$

When specialized to $\mathbf{q}=q$, the affine Hecke algebra arises from functions on the Langlands dual group $\LD{G}$ over a non-Archimedian local field $F$ with residue field $\F_q,$
$$\Hh_{\op{aff},q}\otimes \C= C_c(I\backslash \LD{G}(F)/I, \C).$$
This yields a description of the principal block of smooth $\LD{G}(F)$-representations $\Rep^I(\LD{G}(F))$ in terms of $\Hh_{\op{aff},q}$-modules.  

On the other hand, the generic affine Hecke algebra arises from the geometry of the Springer resolution $\mu:\widetilde{\Nn}\to \Nn\subset \mathfrak{g}^*$ of the nilpotent cone. It is isomorphic to the $G$-theory (=$K$-theory of coherent sheaves) of the Steinberg stack
$$\Hh_{\op{aff},\mathbf{q}} = G_0(\widetilde{\Nn}\times_\Nn\widetilde{\Nn}/(G\times\Gm)).$$

Famously, Kazhdan--Lusztig \cite{kazhdanProofDeligneLanglandsConjecture1987} conceived this geometric description to prove the Deligne--Langlands conjecture which gives a local Langlands correspondence for irreducible smooth representations of $\LD{G}(F)$. 

This motivates the question, which we investigate in this paper, of how to lift this correspondence to a categorical level.

\subsection{Affine Hecke algebra via $K$-motives}
Categories of representations of the \emph{graded} affine Hecke algebra $\Hh_{\op{aff},\mathbf{q}}^{\op{gr}}$ 
can be realized as categories of Springer sheaves, that is, the category generated by $\mu_!(\un)$,
\[\begin{tikzcd}
	{\DMSpr(\Nn/(G\times \Gm),\Q)} & {\DSpr(\Nn/(G\times \Gm),\overline{\Q}_\ell)} \\
	{\Dperf^\Z(\Hh_{\op{aff},\mathbf{q}}^{\op{gr}}\!\otimes\Q)} & {\Dperf(\Hh_{\op{aff},\mathbf{q}}^{\op{gr}}\otimes \overline{\Q}_\ell,d=0).}
	\arrow["\text{\cite{riderFormalityNilpotentCone2013}}"', from=1-2, to=2-2]
    \arrow["\wr", from=1-1, to=2-1]
	\arrow[from=2-1, to=2-2]
	\arrow["\wr", from=1-2, to=2-2]
    \arrow["\text{\cite{eberhardtSpringerMotives2021,eberhardtMotivicSpringerTheory2022}}"', from=1-1, to=2-1]
	\arrow[from=1-1, to=1-2]
\end{tikzcd}\]
A similar description holds for quotients of the affine Hecke algebra at maximal ideals of the center, see \cite{antorFormalityDeligneLanglandsCorrespondence2023}. 
These equivalences reflect the fact that the graded affine Hecke algebra is the Borel--Moore homology of the Steinberg stack, see \cite{lusztigCuspidalLocalSystems1988,lusztigAffineHeckeAlgebras1989}.

To realize representations of the affine Hecke algebra and thereby $\LD{G}(F)$ for all central characters simultaneously, one needs to leave the world of constructible sheaves. 
We propose that \emph{$K$-motives} $\DKbig$, and specifically \emph{reduced $K$-motives} $\DK$, are suitable substitutes for constructible sheaves.\footnote{See \Cref{sec:relationtocoherentspringertheory} for a comparison to a different approach via coherent sheaves}

$K$-motives are a formalism of motivic sheaves, defined for certain stacks, that compute algebraic $K$-theory and admit six functors, thanks to recent work in motivic homotopy theory, for example \cite{cisinskiTriangulatedCategoriesMixed2019a,hoyoisSixOperationsEquivariant2017,hoyoisVanishingTheoremsNegative2019a,khanGeneralizedCohomologyTheories2022}.

Our main result is that the category of representations of the affine Hecke algebra (even with integral coefficients) can be realized in terms of \emph{Springer $K$-motives} on the nilpotent cone.
\begin{theorem*}[\Cref{thm:formalityaffinehecke}] There is an equivalence of categories
$$\Dperf(\Hh_{\op{aff},\mathbf{q}})\simeq \DKSpr(\Nn/(G\times\Gm)).$$
\end{theorem*}
Our methods also generalize the formality results in \cite{eberhardtSpringerMotives2021,eberhardtMotivicSpringerTheory2022} for the graded affine Hecke algebra to integral coefficients.

By specializing at $\mathbf{q}=q,$ we obtain the following categorical version of a local Langlands correspondence for the principal block of $\LD{G}(F)$.
\begin{theorem*}[\Cref{thm:padicrepsviakmotives}] There is an equivalence of categories
    $$\D^b(\Rep^I(\LD{G}(F)))\simeq \DKSpr(\Nn/(G\times\Gm),\C)_{\mathbf{q}=q}.$$
\end{theorem*}
These equivalences further specialize to earlier results on the graded affine Hecke algebra. They are compatible with parabolic induction and the local Langlands parametrization.
\subsection{The formalism of $K$-motives} 
In \Cref{sec:categoriesofkmotives}, we discuss foundational results on equivariant $K$-motives.  
The category of $K$-motives is defined using Hoyois' equivariant motivic stable homotopy category, see \cite{hoyoisSixOperationsEquivariant2017,hoyoisCdhDescentEquivariant2020},
 $$\DKbig(\Xx)=\Mod_{\KGL}(\SH(\Xx))$$
where $\Xx/k$ is a linearly reductive stack, that is, a quotient stacks by the action of a linearly reductive group over a field $k.$

A main distinction between $K$-motives and constructible sheaves is that their mapping spaces compute algebraic $K$-theory and $G$-theory instead of singular co- and Borel--Moore homology.

Furthermore, $K$-motives calculate \emph{genuine equivariant $K$-theory}, in contrast to Borel-style equivariant categories that capture only an `infinitesimal' part, namely the completion at the augmentation ideal of the representation ring. This should be compared to \cite{lusztigAffineHeckeAlgebras1989} where it is shown that the graded affine Hecke algebra is obtained from the affine Hecke algebra via the Atiyah--Segal completion theorem.

For all purposes regarding the affine Hecke algebra, the higher algebraic $K$-theory of the base field $k$ is irrelevant. This is why we use the slightly modified formalism of \emph{reduced $K$-motives}
$\DK(\Xx)$
which, in an appropriate derived sense, are $K$-motives modulo the augmentation ideal $K_{>0}(k).$ This formalism is equipped with six functors using \cite{eberhardtIntegralMotivicSheaves2023}.
\subsection{Formality for reduced $K$-motives}
In \Cref{sec:cellularityandformality}, we prove formality results for $K$-motives arising from `cellular' spaces. We introduce two properties, $(R)$ and $(C)$, for stacks $\Zz.$ They are closely related to the property $(S)$, see \cite{deconciniHomologyZeroSetNilpotent1988}, and, intuitively speaking, ensure that $G_0(\Zz)$ can be computed independently of the higher algebraic $K$-theory of the ground field.

We apply this to the following prototypical situation. Consider a proper map $\mu: \Xx\to \Ss$ between linearly reductive stacks, where $\Xx$ is smooth over $k$. Denote $\Zz=\Xx\times_{\Ss}\Xx.$ Then, we obtain a Hecke-type algebra of $G$-theoretic correspondences
$$\Hh=(G_0(\Zz),\star).$$
Additionally, we define the category of \emph{Springer $K$-motives} 
$$\DKSpr(\Ss)\subset \DK(\Ss)$$ 
as the subcategory generated by $\mu_!(\un)$ and show the following formality result.
\begin{theorem*}[\Cref{thm:propertyrgivesformality}]
    Assume that $\Zz$ fulfills property $(R)$. Then, there is an equivalence
    $$\DKSpr(\Ss)\stackrel{\sim}{\to}\Dperf(\Hh).$$
\end{theorem*}
Similar formality results apply to $\DM$ and constructible sheaves when $\Zz$ satisfies property $(C)$, and we will explain how these formality statements fit together.

In \Cref{sec:kmotivicspringertheory}, we introduce a class of examples for which one can apply the previous theorem. We recall the notion of a \emph{Springer datum} $\mathfrak{S}$ which, put simply, consists of a reductive group with a Borel subgroup $G\supset B$ and a $G$-representation $V\supset U$ with a $B$-stable subspace. To this datum, one associates the map
$$\mu: X=G\times^B U\to S=GU\subset V$$
and the Steinberg variety $Z=X\times_SX.$
We show that $Z/(G\times \Gm)$ satisfies properties $(C)$ and $(R)$. Additionally, we discuss parabolic induction and a Deligne--Langlands classification of irreducibles within the context of Springer data.

Finally, in \Cref{sec:affineheckeandlocallanglands,sec:quiverheckeschuralgebra}, we explain how the affine Hecke algebra and the quiver Hecke/Schur algebra arise from Springer data. We deduce that their representation categories can be realized in terms of Springer $K$-motives.

\subsection{Relation to coherent Springer theory}\label{sec:relationtocoherentspringertheory} In \cite{ben-zviCoherentSpringerTheory2023}, Ben-Zvi, Chen, Helm and Nadler introduce a \emph{coherent Springer theory} (see also \cite{zhuCoherentSheavesStack2021, hellmannDerivedCategoryIwahoriHecke2021} for previous results and conjectures in this direction). They show, for example, that the category of representations of the affine Hecke algebra can be realized in terms of coherent Springer sheaves 
$$\Dperf(\Hh_{\op{aff},\mathbf{q}}\otimes k)\simeq \op{Coh}^{\op{Spr}}(\mathcal{L}(\hat{\Nn}/(G\times\Gm)))$$
on the loop space of the (formal completion) of $\Nn/(G\times\Gm).$ A posteriori, we hence obtain an equivalence
$$\widetilde{\op{ch}}\colon \DKSpr(\Nn/(G\times\Gm),k)\stackrel{\sim}{\to} \op{Coh}^{\op{Spr}}(\mathcal{L}(\hat{\Nn}/(G\times\Gm))).$$

On the level of algebras, the Chern character yields an isomorphism of $K$-theory and the Hochschild homology of the category of coherent sheaves of the Steinberg stack
$$\op{ch}\colon G_0(Z/(G\times\Gm))\otimes k\stackrel{\sim}{\to}HH(\op{Coh}(Z/(G\times\Gm)))$$
as shown in \cite[Corollary 2.26]{ben-zviCoherentSpringerTheory2023}.

On the level of categories, the equivalence $\widetilde{\op{ch}}$ should come from a version of a categorified Chern character between $K$-motives and coherent sheaves on the loop space, introduced in \cite{hoyoisHigherTracesNoncommutative2017}.
The categorified Chern character should enable a direct comparison of similar $K$-motivic and coherent categories in geometric representation theory. We do not investigate this any further in this paper.
\subsection{$K$-motives and Langlands duality}
The current paper is part of an ongoing program to construct equivalences between $K$-motivic and constructible categories on Langlands dual objects. We outline some other results and conjectures in this direction.

In \cite{eberhardtKMotivesKoszulDuality2019}, we proposed a $K$-motivic approach to classical Koszul duality for Langlands dual flag varieties, see \cite{soergelKategoriePerverseGarben1990,beilinsonKoszulDualityPatterns1996}. By combining Soergel-theoretic descriptions for monodromic sheaves, see \cite{taylorUncompletingSoergelEndomorphismensatz2023}, and $K$-motives, see \cite{eberhardtKtheorySoergelBimodules2022}, one obtains the following global and ungraded form of Koszul duality for finite flag varieties.
\begin{theorem*}[Global Koszul duality]
    There is an equivalence
    $$\DKbig(B\backslash G/ B, \Q)\simeq \D_{\op{mon}}(\LD{U}\backslash \LD{G}(\C) /\LD{U}, \Q)$$
    between $K$-motives on the finite Hecke stack and monodromic sheaves on the Langlands dual. The equivalence exchanges intersection $K$-motives and free monodromic tilting sheaves.
\end{theorem*}
The equivalence is linear over the representation ring $R$ of the maximal torus $T\subset G.$
By completing at the augmentation ideal $I\subset R$, this specializes to the Koszul duality of \cite{bezrukavnikovKoszulDualityKacMoody2013} involving equivariant and unipotent monodromic sheaves.

Soergel's conjecture \cite{soergelLanglandsPhilosophyKoszul2001} 
predicts a similar Koszul duality for real reductive groups.
Denote by $\theta$ the holomorphic Cartan involution for a real form of a complex reductive group $G.$ 
By passing to homotopy fixed points of $\theta$, one obtains the following version of the conjecture.
\begin{conjecture*}[Global Soergel's conjecture]
There is an equivalence 
    $$\DKbig([B\backslash  G/B]^{h\theta})\simeq \op{D}_{\op{mon}}([\LD{U}\backslash \LD{G}/\LD{U}]^{h\LD{\theta}}(\C))$$
    between $K$-motives and monodromic sheaves on homotopy fixed points of corresponding Galois actions on the Hecke stacks.
\end{conjecture*}
In \cite{ben-zviCoherentConstructibleLocal2023} an overview on similar equivalences involving coherent sheaves on loop spaces instead of $K$-motives is given. As outlined in \Cref{sec:relationtocoherentspringertheory}, these approaches should be related by a categorical Chern character.

Another direction is the quantum $K$-theoretic Satake, relating equivariant $K$-motives on the affine Grassmannian with representations of the Langlands dual quantum group.
\begin{conjecture*}[Quantum geometric Satake]
    There is an equivalence
    $$\DK((G\times \Gm)\backslash\!\operatorname{Gr}_G)\stackrel{\sim}{\to}\D^{b}_{U_q(\LD{\mathfrak{g}})}(\mathcal{O}_q(\LD{G})),$$
\end{conjecture*}
For $\op{SL}_n$, a similar conjecture was proven in \cite{cautisQuantumKtheoreticGeometric2018}.

\subsection{Acknowledgements} We thank Harrison Chen, Arnaud Eteve, Quoc Ho and Marc Hoyois for helpful discussions. Moreover, we thank Jonas Antor for valuable input regarding the proof of the formality result.
The author was supported by Deutsche Forschungsgemeinschaft (DFG), project number 45744154, Equivariant K-motives and Koszul duality.
\section{Categories of $K$-motives on stacks}\label{sec:categoriesofkmotives}
This section introduces the category of (reduced) $K$-motives on quotient stacks and discusses its functorial properties and relationships with other categories of equivariant motives and sheaves.
\subsection{Quotient stacks}
Fix the base scheme $\point=\Spec(k)$ for an algebraically closed field $k$ of characteristic zero. We call a stack $\Xx/\point$ \emph{\nice} if it is isomorphic to a quotient stack $X/G$ where
\begin{enumerate}
    \item $G$ is a linearly reductive group,
    \item $X$ is $G$-quasi-projective, that is, admits an embedding in a projectivized representation of $G$, and of finite type.
\end{enumerate}
We write $\Bc G=\point/G$ for the classifying stack.
\subsection{$K$-motives on \nice\ stacks}
For a \nice\ stack $\Xx$, we denote by  $\SH(\Xx)$ the stable motivic homotopy category associated to $\Xx$ as defined in \cite{hoyoisSixOperationsEquivariant2017}. 

By \cite{hoyoisCdhDescentEquivariant2020}, there is a ring spectrum $\KGL_\Xx$ in $\SH(\Xx)$ representing $\A^1$-homotopy invariant algebraic $K$-theory. The category of \emph{$K$-motives} on $\Xx$ is defined by 
$$\DKbig(\Xx)=\Mod_{\KGL_\Xx}(\SH(\Xx)).$$
In \Cref{sec:reducedkmotives} we will define a reduced variant $\DK$ of $\DKbig$. In \cite{hoyoisSixOperationsEquivariant2017,khanGeneralizedCohomologyTheories2022} it is shown that the system of categories $\Xx\mapsto\DKbig(\Xx)$ admits
\begin{enumerate}
    \item $f^*,f_*,$ for any morphism,
    \item $f_!,f^!$ for representable morphisms and
    \item bifunctors $\iHom,\otimes$
\end{enumerate}
which fulfill the axioms of a six functor formalism, such as base change, localization and projection formulae. 

\begin{remark}A special feature of $K$-motives is \emph{Bott periodicity}: The functor $(1)[2]$ is isomorphic to the identity.
Here, $(1)[2]$ denotes the cofiber of 
$\id\to \pi_*\pi^*$ for $\pi: \Proj^1_\Xx\to \Xx.$ This is reflected in the fact that Verdier duality comes as an equivalence of functors $f^!\simeq f^*$ for $f$ smooth representable. 
\end{remark}

The category of $K$-motives computes algebraic $K$-theory and $G$-theory in the following sense.
For $\Xx/\point$ a smooth \nice\ stack there is an equivalence of ring spectra
\begin{align}\label{eq:KtheoryasmappingspaceinDKbig}
    \Maps_{\DKbig(\Xx)}(\un,\un) &\simeq K(\Xx),\\
\intertext{where $K(\Xx)$ denotes $K$-theory spectrum of the category of perfect complexes on $\Xx$, see \cite[below Definition 5.1]{hoyoisCdhDescentEquivariant2020}.
If $\Xx$ is not necessarily smooth and $f:\Xx\to \Yy$ is a representable map to a regular \nice\ stack $\Yy$, then there is an equivalence of spectra
}
	\label{eq:GtheoryasmappingspaceinDKbig}
    \Maps_{\DKbig(\Xx)}(\un,f^!\un)& \simeq G(\Xx),
\end{align}
where $G(\Xx)$ denotes the $K$-theory spectrum of the category of coherent sheaves on $\Xx$, see \cite[Remark 5.7]{hoyoisCdhDescentEquivariant2020}.
\subsection{Reduced $K$-motives} \label{sec:reducedkmotives} 
For the applications we have in mind in this paper, we are not interested in the higher algebraic $K$-theory of the base field $k.$ We will hence `remove' these groups using the reduction methods from \cite{eberhardtIntegralMotivicSheaves2023}. There, a similar procedure is performed in the context of $\DM$. We briefly sketch the construction and refer to \emph{loc.\ cit.}\ for details and references.

Recall that $\point=\Spec(k)$ where $k$ is an algebraically closed field of characteristic zero and $K(\point)$ is connective, that is, $K_i(\point)=\pi_{i}(K(\point))=0$ for $i<0$.

We consider the stable idempotent closed  subcategory $\DKbigconst(\point)\subset \DKbig(\point)$ generated by the unit object $\un$. 
Since $K(\point)$ is connective, $\DKbigconst(\point)$ admits a weight structure $w$ whose heart $\DKbigconst(S)^{w=0}$ is the additive idempotent closed subcategory generated by $\un$, see \cite[Theorem 4.3.2]{bondarkoWeightStructuresVs2010}. The homotopy category of $\DKbigconst(\point)^{w=0}$ is equivalent to the additive category of finitely generated projective modules of $K_0(\point)=\Z$
$$\ho\DKbigconst(\point)^{w=0}\simeq \Proje_{\op{f.g.}}(K_0(\point)).$$ 

The weight structure induces a weight complex functor
$$\DKbigconst(\point)\to \Ch^b(\ho\DKbigconst(\point)^{w=0})\simeq \Dperf(K_0(\point))$$
which comes from pushforward along the truncation map $K(\point)\to K_0(\point).$ The functor is symmetric monoidal by \cite{aokiWeightComplexFunctor2020}.
For \nice\ stacks $f:\Xx\to\point$ we define the category of \emph{reduced $K$-motives} on $\Xx$ via the Lurie tensor product
$$\DK(\Xx)=\DKbig(\Xx)\otimes_{\DKbigconst(\point)}\Dperf(K_0(\point)).$$
Here, we $N\in \DKbigconst(\point)$ acts on $M\in\DKbig(\Xx)$ via $M\otimes f^*N.$
We denote by
\begin{align*}
	K(\Xx)_\red=K(\Xx)\otimes_{K(\point)}K_0(\point)\text{ and }
	G(\Xx)_\red=G(\Xx)\otimes_{K(\point)}K_0(\point)
\end{align*}
the \emph{reduced $K$-theory} and \emph{reduced $G$-theory} of $\Xx$ and by $K_i(\Xx)_\red$ and $G_i(\Xx)_\red$ the $i$-th homotopy group, respectively. Using \eqref{eq:KtheoryasmappingspaceinDKbig} and \eqref{eq:GtheoryasmappingspaceinDKbig} and the properties of the Lurie tensor product, there are equivalences
\begin{align}
    \Maps_{\DK(\Xx)}(\un,\un) &\simeq K(\Xx)_\red\text{ and}\\
    \Maps_{\DK(\Xx)}(\un,f^!\un)& \simeq G(\Xx)_\red\label{eq:GtheoryasmappingspaceinDK}
\end{align}
where in the first equivalence we assume that $\Xx/\point$ is smooth and in the second that $f:\Xx\to \Yy$ is a representable map to a regular \nice\ stack $\Yy.$

In \cite{eberhardtIntegralMotivicSheaves2023} (in the context of $\DM$, however the proofs translate verbatim to $\DKbig$) it is shown that $\Xx\mapsto\DK(\Xx)$ inherits a six functor formalism compatible with the reduction functor $\DKbig(\Xx)\to \DK(\Xx).$ 
\subsection{Realization and completion}\label{sec:realisation}
The category of (reduced) $K$-motives can be compared to Borel-style equivariant motives and sheaves via the following diagram
\[\begin{tikzcd}[column sep=small, row sep = tiny]
	{\D(\Xx_{\Borel},\Lambda)} & {\DM_{(\red)}(\Xx_{\Borel},{\Q})} & {\DKbig_{(\red)}(\Xx_{\Borel},\Q)} & {\DKbig_{(\red)}(\Xx)} \\
	{\displaystyle\prod_nH^{2n}(\Xx_{\Borel},\Lambda)} & {\displaystyle\prod_n\CH^n(\Xx_{\Borel},\Q)_{(\red)}} & {K_0(\Xx_{\Borel})_{(\red)}\otimes \Q} & {K_0(\Xx)^\wedge_{(\red),I}\otimes \Q}
	\arrow["\pi"',from=1-4, to=1-3]
	\arrow["\rho"', from=1-2, to=1-1]
	\arrow["v", from=1-2, to=1-3]
	\arrow[from=2-2, to=2-1]
    \arrow["\sim"', from=2-3, to=2-2]
    \arrow["\sim"', from=2-4, to=2-3]
\end{tikzcd}\]
which, for smooth \nice\ stacks $\Xx/\point$, on the level of cohomology induces, from left to right, the cycle class, Chern character and completion map. Let us explain the notation.
\begin{enumerate}
    \item 
	
	For a six functor formalism $\DD$, defined for schemes, we denote by $$\DD(\Xx_{\Borel})=\lim_{(U,u)\in \operatorname{Lis}_\Xx}\DD(U)$$ the \emph{lisse extension} defined by the limit along $*$-pullbacks over the category $\operatorname{Lis}_\Xx$ of all smooth morphisms $u:U\to \Xx$ where $U$ is a qcqs algebraic space, see \cite{khanGeneralizedCohomologyTheories2022}. For most practical matters and the applications discussed here, this coincides with similar definitions of Borel-style equivariant sheaves, see for example \cite{bernsteinEquivariantSheavesFunctors1994, soergelEquivariantMotivesGeometric2018, richarzIntersectionMotiveModuli2020}.
    \item By $\DKbig_{(\red)}$ and $\DM_{(r)}$ we denote either the categories of  ($K$-)motives and reduced ($K$-)motives, see \Cref{sec:reducedkmotives} and \cite{eberhardtIntegralMotivicSheaves2023}.
    \item By $\D$ we denote either the derived category of $\ell$-adic or Betti sheaves, that is, sheaves on the complex points equipped with the metric topology, with $\Lambda=\overline\Q_\ell$ and $\Lambda=\C,$ respectively.
    \item The \emph{completion functor} $\pi$ is induced by the functor $u^*:\DKbig_{(\red)}(\Xx)\to \DKbig_{(\red)}(U)$ for $(U,u)\in\operatorname{Lis}_\Xx.$
    \item The \emph{Beilinson realization functor} $v$ exists with rational coefficients and induces the Chern character map from $K$-theory to Chow groups, see \cite{riouAlgebraicKtheoryA1homotopy2010,cisinskiTriangulatedCategoriesMixed2019a,eberhardtKMotivesKoszulDuality2019}.
    \item The \emph{$\ell$-adic} or \emph{Betti realization functor} $\rho$ induces the cycle class map from Chow groups to $\ell$-adic/Betti cohomology, see \cite{cisinskiEtaleMotives2016,ayoubNoteOperationsGrothendieck2009} as well as \cite[Section 5]{eberhardtIntegralMotivicSheaves2023} for reduced motives.
    \item By $K_0(\Xx_{\Borel})_{(\red)}$, $\CH^n(\Xx_{\Borel})_{(\red)}$ and $H^{2n}(\Xx_{\Borel})$ we denote the Borel-style equivariant cohomology constructed using algebraic approximations of the classifying space $\Xx_{\Borel}=\lim_i U_i\times^G\!X$ in the sense of Totaro, \cite{totaroChowRingClassifying1999}.
\end{enumerate}
\begin{remark}
On the level of cohomology, the map from the  Bredon-style equivariant $K$-theory to the Borel-style equivariant $K$-theory induces an isomorphism of the former completed at the augmentation ideal $I=\ker(K_0(\Bc G)\to K_0(\point))$ with the latter via the Atiyah--Segal completion theorem, \cite{atiyahEquivariantTheoryCompletion1969,thomasonEquivariantAlgebraicVs1988,totaroChowRingClassifying1999}. 
Hence, $\DKbig_{(\red)}(\Xx_{\Borel})$ can be seen as a \emph{completed} and \emph{infinitesimal} version of $\DKbig_{(\red)}(\Xx).$ 

Moreover, we interpret $\DM_{(\red)}(\Xx_{\Borel},\Q)$ as a graded version of $\DKbig_{(\red)}(\Xx_{\Borel},\Q)$ where $(1)[2]$ shifts and $v$ forgets the grading. Similarly, $\DM_{(\red)}(\Xx_{\Borel},\Q)$ admits an additional grading with respect to $\Dc(\Xx_{\Borel},\Lambda)$ where $(1)$ shifts and $\rho$ forgets the grading.

This motivates the slogan that $K$-motives are a global version of equivariant sheaves. In particular, one may consider completions $\DK(\Xx)^\wedge_\mathfrak{m}$ at different maximal ideals $\mathfrak{m}\subset K_0(\Bc G)$ which are $(1)[2]$-periodic $\mathfrak{m}$-twisted equivariant mixed sheaves.
\end{remark}
\section{Cellularity and formality}\label{sec:cellularityandformality}
We introduce cellularity conditions for stacks and how they imply formality for Springer $K$-motives.
\subsection{Properties $(R)$ and $(C)$} We introduce properties of stacks which ensure that reduced $G$-theory $G(\Zz)_{\red}$ is equivalent to $G_0(\Zz).$
\begin{definition} Let $\Zz$ be a stack. Denote by $G(\Zz)\to G(\Zz)_{\red}$ the natural reduction map.  We say that $\Zz$ has property $(R)$ if
	\begin{enumerate}
        \item[$(R1)$] $G_i(\Zz)_{\red}=0$ for $i\neq0$ and
        \item[$(R2)$] the map $G_0(\Zz)\to G_0(\Zz)_{\red}$ is an isomorphism.
    \end{enumerate}
\end{definition}
The fundamental example for the property $(R)$ comes from classifying stacks.
\begin{proposition}\label{prop:propertyRforBG}
	Let $G$ be reductive and $U$ a unipotent linear algebraic group over $k$. Then $\Bc(G\rtimes U)$ has property $(R).$
\end{proposition}
\begin{proof}
	There are equivalences $G(\Bc (G\rtimes U))\stackrel{\sim}{\rightarrow}K(\Bc(G\rtimes U))\stackrel{\sim}{\to}K(\Bc G)$ of $K(\point)$-modules by \cite[Lemma 4.6]{thomasonLefschetzRiemannRochTheoremCoherent1986} and \cite[Theorem 5.7]{thomasonXXAlgebraicKTheory2016}.
	Since $G$ is reductive and $k$ an algebraically closed field of characteristic $0$, $\Rep(G)$ is a direct sum of categories equivalent to finite-dimensional $k$-vector spaces. Since $K$-theory commutes with colimits, this yields an equivalence of $K(\point)$-modules $K_0(\Bc G) \otimes_{K_0(\point)}K(\point)\stackrel{\sim}{\to} K(\Bc G).$ Hence, $G(\Bc(G\rtimes U))_{\red}=G_0(\Bc G).$ This implies the statement.
\end{proof}

The following allows to inductively deduce property $(R)$ using open and closed decompositions and vector bundles.
\begin{lemma}\label{lem:propertyRinductive} Let $\Zz$ be a stack.
	\begin{enumerate}
		\item If $\Vv\to \Zz$ is a vector bundle and $\Zz$ has property $(R)$, then so does $\Vv.$
		\item If $\Uu\subset \Zz$ is an open substack and $\Uu$ and $\Zz-\Uu$ have property $(R)$, then so does $\Zz.$
	\end{enumerate}
\end{lemma}
\begin{proof}
(1) follows from \Cref{prop:propertyRforBG} and homotopy invariance, see \cite[Theorem 4.1]{thomasonXXAlgebraicKTheory2016}. 
For (2), assume that $\Uu$ and $\Zz-\Uu$ have property $(R).$ Consider the map of long exact sequences
\[\begin{tikzcd}[column sep=small]
	{G_1(\Uu)} & {G_0(\Zz-\Uu)} & {G_0(\Zz)} & {G_0(\Uu)} & {G_{-1}(\Zz-\Uu)} \\
	{G_1(\Uu)_{\red}} & {G_0(\Zz-\Uu)_{\red}} & {G_0(\Zz)_{\red}} & {G_0(\Uu)_{\red}} & {G_{-1}(\Zz-\Uu)_{\red}}
	\arrow[from=2-2, to=2-3]
	\arrow[from=2-3, to=2-4]
	\arrow[from=2-4, to=2-5]
	\arrow[from=2-1, to=2-2]
	\arrow[from=1-3, to=2-3]
	\arrow[from=1-2, to=1-3]
	\arrow["\wr", from=1-2, to=2-2]
	\arrow[from=1-3, to=1-4]
	\arrow[from=1-4, to=1-5]
	\arrow["\wr", from=1-4, to=2-4]
	\arrow[from=1-5, to=2-5]
	\arrow[from=1-1, to=2-1]
	\arrow[from=1-1, to=1-2]
\end{tikzcd}\]
where the second map and fourth vertical maps are isomorphisms and $G_1(\Uu)_{\red}$ is zero by assumption. Moreover, $G_{-1}(\Zz-\Uu),$ is zero by the definition of $G$-theory. By the five lemma, the third vertical map is an isomorphism, so $\Zz$ fulfills $(R2)$. Similarly, $(R1)$ follows from the assumption and the bottom long exact sequence.
\end{proof}
We will deduce property $(R)$ from the following cellularity property for stacks.
\begin{definition}
	Let $\Zz$ be a \nice\ stack. We say that $\Zz$ has property $(C)$ if it admits a filtration $\Zz=\Zz^n\supset \Zz^{n-1}\supset\cdots\Zz^1\supset \Zz^0=\emptyset$ over $\Ss$ such that $\Zz^i\to \Zz$ is a closed immersion, and each $\Vv^i=\Zz^i-\Zz^{i-1}$ is a vector bundle over $\Bc (G_i\ltimes U_i)$ for $G_i$ reductive and $U_i$ unipotent.
\end{definition}
\begin{corollary}\label{cor:cimpliesr}
	Property $(C)$ implies property $(R).$
\end{corollary}
\begin{proof}
	Combine \Cref{prop:propertyRforBG} and \Cref{lem:propertyRinductive}.
\end{proof}
\subsection{Springer motives and correspondences}\label{sec:springercorrespondences}
Let $\mu:\Xx\to \Ss$ be a proper map of \nice\ stacks over $\point$ such that $\Xx/k$ is smooth. Denote $\Zz=\Xx\times_\Ss\Xx.$
\begin{definition}
	The stable and idempotent closed subcategory generated by the object $\mu_!\un_\Xx$
$$\DKbigSpr_{(\red)}(\Ss)=\langle \mu_!\un_\Xx\rangle_{stb,\inplus} \subset \DKbig_{(\red)}(\Ss)$$
is called the category of (reduced) \emph{Springer $K$-motives.}
\end{definition}
The category of Springer $K$-motives is hence governed by the endomorphism mapping space of $\mu_!\un_\Xx,$ which can be computed as follows.
\begin{proposition}\label{prop:mappingsviagtheoryofsteinberg}
    There is an equivalence of spectra
\begin{align*}
	\Maps_{\DKbig_{(\red)}(\Ss)}(\mu_!\un_\Xx,\mu_!\un_\Xx)&\simeq G(\Zz)_{(\red)}
\end{align*}
\end{proposition}
\begin{proof} 
	Using base change for the Cartesian diagram
	\[\begin{tikzcd}
		\Zz & \Xx \\
		\Xx & \Ss
		\arrow["{\pi_1}"', from=1-1, to=2-1]
		\arrow["{\pi_2}", from=1-1, to=1-2]
		\arrow["p", from=1-2, to=2-2]
		\arrow["p", from=2-1, to=2-2]
	\end{tikzcd}\]
as well as the fact that $\Xx$ is regular, from \eqref{eq:GtheoryasmappingspaceinDKbig} and \eqref{eq:GtheoryasmappingspaceinDK} we obtain equivalences of spectra
\begin{align*}
	\Maps_{\DKbig_{(\red)}(\Ss)}(p_!\un_\Xx,p_!\un_\Xx)&\simeq\Maps_{\DKbig_{(\red)}(\Zz)}(\pi_1^*\un_\Xx,\pi_2^!\un_\Xx)\\
	&\simeq\Maps_{\DKbig_{(\red)}(\Zz)}(\un_\Zz,\omega_{\Zz/\Xx})\\
	&\simeq G(\Zz)_{(\red)}.\qedhere
\end{align*}
\end{proof}
On the level of homotopy groups, the above equivalence also respects the algebra structure, as the following statement shows.
\begin{proposition}\label{prop:compositionequalsconvolution} There are isomorphisms of algebras
    \begin{align*}\End_{\DKbig_{(\red)}(\Ss)}(\mu_!\un_\Xx)&\cong (G_0(\Zz)_{(\red)},\star).
	\end{align*}
    The right-hand side is equipped with the product $M\star N=p_*\Delta^!(M\boxtimes N)$ where $\Delta^!$ is the refined Gysin morphism and $p_*$ for the diagram
\[\begin{tikzcd}[column sep=small, row sep=0]
	{\Xx\times_{\Ss}\Xx\times \Xx\times_{\Ss}\Xx} && {\Xx\times_{\Ss}\Xx\times_{\Ss}\Xx} && {\Xx\times_{\Ss}\Xx} \\
	& \Delta \\
	{\Xx\times\Xx\times \Xx\times\Xx} && \Xx\times\Xx\times\Xx.
	\arrow[from=1-3, to=1-1]
	\arrow["p", from=1-3, to=1-5]
	\arrow[from=1-1, to=3-1]
	\arrow[from=1-3, to=3-3]
	\arrow[from=3-3, to=3-1]
\end{tikzcd}\]\end{proposition}
\begin{proof}This is shown in \cite{fangzhouBorelMooreMotivic2016} for $\DM$. The proof works the same here.
\end{proof}
\subsection{Property $(R)$ and formality of Springer $K$-motives}
We denote by $\Hh=(G_0(\Zz),\star)$ the algebra of $G$-theoretic correspondences associated to $\Zz,$ see \Cref{prop:compositionequalsconvolution}. This is an ordinary algebra with no higher structure.

\begin{theorem}
	\label{thm:propertyrgivesformality}
    If $\Zz$ fulfills property $(R)$, then there is an equivalence
    $$\DKSpr(\Ss)\to \Dperf(\Hh).$$
\end{theorem}
\begin{proof} \Cref{prop:mappingsviagtheoryofsteinberg} and property $(R)$ show that
    $$\Maps_{\DK(\Yy)}(\mu_!\un_\Xx,\mu_!\un_\Xx)\simeq G_0(\Zz)$$
    is concentrated in degree zero. Hence, there is a weight structure on $\DKSpr(\Ss)$ with heart generated by finite direct sums and retracts of $p_!\un_\Xx.$ Moreover, the weight complex functor yields an equivalence
    $$\DKSpr(\Ss)\to \Ch^b(\ho\DKSpr(\Ss)^{w=0})\cong \Dperf(\End_{\DK(\Ss)}(\mu_!\un_\Xx)).$$
    Now, \cref{prop:compositionequalsconvolution} computes $\End_{\DK(\Ss)}(\mu_!\un_\Xx)=(G_0(\Zz),\star)=\Hh.$
\end{proof}

\subsection{Realization and completion}\label{sec:realandcompletionspringermotives} Assume that $\Zz$ fulfills property $(C).$ Then, \Cref{thm:propertyrgivesformality} applies verbatim to the for Borel-style reduced $K$-motives and motives, see  \Cref{sec:realisation}.
Hence, for example, there is an equivalence
$$\DMSpr_{\red}(\Ss_{\Borel})\stackrel{\sim}{\to} \Dperf^\Z(\Hh^{\op{gr}})$$
where $\Hh^{\op{gr}}=(\bigoplus_n \CH^n(\Zz_{\Borel}),\star)$ arises from the equivariant Chow groups and
$$\DMSpr_{\red}(\Ss_{\Borel})=\genbuild{\mu_!\un(n)[2n]}{n\in \Z}_{stb,\inplus}\subset \DM_{\red}(\Ss_{\Borel})$$
is the category of Borel-style equivariant Springer motives.
Similar result with rational coefficients were obtained in \cite{eberhardtSpringerMotives2021, eberhardtMotivicSpringerTheory2022}.
For $\ell$-adic or Betti sheaves and $\Lambda=\overline\Q_\ell$ or $\Lambda=\C$, respectively, one may deduce along the lines of \cite{riderFormalityNilpotentCone2013,antorFormalityDeligneLanglandsCorrespondence2023} that there is an equivalence for $\ell$-adic or Betti Springer sheaves
$$\DSpr(\Ss,\Lambda)\stackrel{\sim}{\to} \Dperf(\Hh^{\op{gr}}\otimes\Lambda,d=0)$$
where one interprets $\Hh^{\op{gr}}_\Lambda$ as $dg$-algebra with trivial differential. In total, we can compare the various formality results for Springer sheaves/motives via the diagram
\[\begin{tikzcd}[column sep=14pt]
	{\DSpr(\Ss_{\Borel},\Lambda)} & {\DMSpr_{\red}(\Ss_{\Borel},\Q)} & {\DKSpr(\Ss_{\Borel},\Q)} & {\DKSpr(\Ss)} \\
	{\Dperf(\Hh^{\op{gr}}\otimes\Lambda,d=0)} & {\Dperf^\Z(\Hh^{\operatorname{gr}}\otimes\Q)} & {\Dperf(\Hh_{I}^\wedge\otimes\Q)} & {\Dperf(\Hh)}
	\arrow["\pi"', from=1-4, to=1-3]
	\arrow["v", from=1-2, to=1-3]
	\arrow["\wr", from=1-2, to=2-2]
	\arrow[from=2-4, to=2-3]
	\arrow["\wr", from=1-4, to=2-4]
	\arrow["\wr", from=1-3, to=2-3]
	\arrow["v", from=2-2, to=2-3]
	\arrow["\rho"', from=1-2, to=1-1]
	\arrow["\wr", from=1-1, to=2-1]
	\arrow[from=2-2, to=2-1]
\end{tikzcd}\]
The diagram involves the Hecke algebras
\[\begin{tikzcd}[column sep=9.5pt]
	{\Hh^{\operatorname{gr}}\otimes\Lambda} & {\Hh^{\operatorname{gr}}\otimes\Q} & {\overline{\Hh}^{\operatorname{gr}}\otimes\Q} & {\Hh^\wedge_{I}\otimes \Q} & \Hh \\
	{\bigoplus\limits_n\!H^{2n}(\Zz_{\Borel},\Lambda)} & {\bigoplus\limits_n\CH^n(\Zz_{\Borel},\Q)} & {\prod\limits_n\CH^n(\Zz_{\Borel},\Q)} & {G_0(\Zz_{\Borel})\!\otimes\!\Q} & {G_0(\Zz)}
	\arrow["\pi"', from=2-5, to=2-4]
	\arrow[shift right, Rightarrow, no head, from=1-2, to=2-2]
	\arrow[hook, from=2-2, to=2-3]
	\arrow[hook, from=1-2, to=1-3]
	\arrow["\sim"', from=1-4, to=1-3]
	\arrow[from=1-5, to=1-4]
	\arrow[Rightarrow, no head, from=1-5, to=2-5]
	\arrow[Rightarrow, no head, from=1-4, to=2-4]
	\arrow[Rightarrow, no head, from=1-3, to=2-3]
	\arrow["\RR"', from=2-4, to=2-3]
	\arrow[from=2-2, to=2-1]
	\arrow[from=1-2, to=1-1]
	\arrow[Rightarrow, no head, from=1-1, to=2-1]
\end{tikzcd}\]
where $H^{2n}(\Zz_{\Borel})_\Lambda$, $\CH^n(\Zz_{\Borel})$ and $G_0(\Zz_{\Borel})$ we denote the Borel-style equivariant groups, see \Cref{sec:realisation}. Moreover, $\operatorname{RR}$ denotes the Riemann--Roch map, see \cite[Theorem 5.11.11]{chrissRepresentationTheoryComplex2010}, which is the product of the Chern character map by the inverse of the Todd class of $\Xx.$ This twist by the Todd class is necessary in order to commute with convolution. 
The Todd class naturally arises from the comparison of orientations of the spectra in $\SH$ representing motivic cohomology and $K$-theory, see \cite{degliseBivariantTheoriesMotivic2018,degliseOrientationTheoryArithmetic2019}.
Since $\Zz$ fulfills property $(C)$, after tensoring with $\Lambda$ the cycle class map becomes an isomorphism as well, as a straightforward induction argument shows.

\section{$K$-motivic Springer theory}\label{sec:kmotivicspringertheory}
We discuss a general framework for $K$-motivic Springer theory, which encompasses formality, parabolic induction, central characters and parametrizations of irreducibles encompassing the affine Hecke algebra and similar algebras.
\subsection{Springer datum} \label{sec:springerdatum}
Let $G\supset B\supset T$ be a reductive group with a Borel subgroup and maximal torus. Let $W=N_G(T)/T$ be the Weyl group. 
\begin{definition}
	A \emph{Springer datum} $\mathfrak{S}=(G,\{P_i\},V,\{V_i\},I,G')$ consists of
\begin{enumerate}
    \item a representation $V$ of a reductive group $G$ and a finite index set $I,$
    \item for each $i\in I$ a subspace $V_i\subset V$ and parabolic subgroup $B\subset P_i\subset G$ such that $V_i$ is $P_i$-stable,
    \item $G'=G\times D$ or $\{1\}\times D$ for a subgroup $D\subset \Gm.$
\end{enumerate}
\end{definition}
\begin{remark}
	When $I$ is a singleton set, we simply write $\mathfrak{S}=(G,P,V,U,G')$ for $U\subset V$ and $P\subset G.$
\end{remark}

To this datum we associate the diagram
\[\begin{tikzcd}
	\bigsqcup\limits_{i\in I}G/P_i & {X=\bigsqcup\limits_{i\in I}X^i} & {S=\im(\mu)\subset V}
	\arrow["\mu", from=1-2, to=1-3]
	\arrow["\pi"', from=1-2, to=1-1]
\end{tikzcd}\]
where $X^i=G\times^{P_i} V_i$, $\pi$ and $\mu$ map $[g,v]$ to $gP_i$ and $gv,$ respectively. 
There is a natural action of $G\times \Gm$ on the diagram, where $t\in\Gm$ acts via dilation by $t^{-1}$ on $V.$

The \emph{Steinberg variety} is defined as
$Z=X\times_SX.$ We denote the irreducible components by $Z^{i,j}=X^i\times_{S}X^{j}$ for $i,j\in I.$
We obtain a natural map $Z^{i,j}\to G/P_i\times G/P_{j}$ and denote for $w\in W$ the preimage of the Bruhat cell $Y_w^{i,j}=G(eP_i,wP_{j})$ by $Z^{i,j}_w.$
The action of $G\times \Gm$ on $Z$ leaves the subspaces $Z^{i,j}_w$ invariant.
\begin{lemma}\label{lem:springerdatumgivescellularfibration}
    \begin{enumerate}[wide, labelwidth=!, labelindent=0pt]
		\item The map $\mu$ is proper and $X^i/\point$ is smooth.
		\item The stack $Z/G'$ fulfills property $(C).$
	\end{enumerate}
\end{lemma}
\begin{proof} (1) is standard. 
For (2) it is enough to show that each $Z_w^{i,j}/H$ fulfills property $(C).$
Now, $Z_w^{i,j}\to Y_w^{i,j}$ is a vector bundle with fiber $V_i\cap w(V_j).$ So we reduced the statement to showing that $Y^{i,j}_w/G'$ fulfills property $(C)$. 

If $G'=G\times D$ we obtain $Y_w^{i,j}/G'=\Bc(P_i\cap P_j^w\times D)$ and $P_i\cap P_j^w\times D$ can be written as $H\ltimes U$ for $H$ reductive and $U$ unipotent and the statement follows. 

If $G'=\{1\}\times D,$ $Y_w^{i,j}/G'=Y_w^{i,j}\times \Bc D.$ Then the projection $Y_w^{i,j}\to G/P_i$ is a $G$-equivariant bundle with fiber $P_iwP_j/P_j.$ The product of the Bruhat decompositions of $G/P_i$ and $P_iwP_j/P_j$ yields an affine stratification of $Y_w^{i,j}.$ The statement follows.
\end{proof}
\begin{remark} We use this general setup here since it also applies, for example, to quiver Hecke algebras (also called KLR algebras) and quiver Schur algebras. We refer to  \cite{sauterSurveySpringerTheory2013,sauterCellDecompositionsQuiver2016} for an overview on `Springer theories'.  
\end{remark}
\subsection{Formality for a Springer datum} 
For a Springer datum $\mathfrak{S}$ we consider the category of Springer motives, as defined in \Cref{sec:springercorrespondences},  associated to the map $\mu:X/G'\to S/G'$
$$\DKSpr(\mathfrak{S})=\DKSpr(S/G')=\langle \mu_!\un_{X/G'}\rangle_{stb,\inplus} \subset \DK(S/G').$$ 
Moreover, we denote by $\Hh(\mathfrak{S})=(G_0(Z/G'),\star)$ the associated \emph{Hecke algebra}.
\begin{theorem}\label{thm:steinbergformality} There is an equivalence of categories
    $$\DKSpr(\mathfrak{S})\stackrel{\sim}{\to}\Dperf(\Hh(\mathfrak{S})).$$
\end{theorem}
\begin{proof} 
	By \Cref{lem:springerdatumgivescellularfibration}, $Z/G'$ fulfills Property $(C)$ and hence, by \Cref{cor:cimpliesr}, Property $(R).$ So \Cref{thm:propertyrgivesformality} implies the statement.
\end{proof}
\subsection{Parabolic Induction}\label{sec:parabolicinduction} We explain how this relates to parabolic induction. For simplicity, we fix a Springer datum $\mathfrak{S}=(G,B,V,U,G')$ for $G'=G\times D$. Then $X=G\times^BU$ and $S=\im(\mu)=GU\subset V.$

Let $B\subset P\subset G$ be a standard parabolic and $W\subset U$ is an $R_u(P)$-stable subspace, where $R_u(P)\subset P$ is the unipotent radical of $P.$ Denote by $L=P/R_u(P),$ $B_L=B/R_u(P)$ and $V_L=V/W$, $U_L=U/W$ as well as $P'=P\times D,$ $L'=L\times D$ and so on.

This yields a Springer datum $\mathfrak{S}_L$ for $L$ and a commutative diagram
\[\begin{tikzcd}
	{X=G\times^BU} & {X_P=P\times^BU} & {X_L=L\times^{B_L}U_L} \\
	S & {S_P} & {S_L.}
	\arrow["{i'}"', from=1-2, to=1-1]
	\arrow["\mu", from=1-1, to=2-1]
	\arrow["{\mu_P}", from=1-2, to=2-2]
	\arrow["i"', from=2-2, to=2-1]
	\arrow["p", from=2-2, to=2-3]
	\arrow["{p'}", from=1-2, to=1-3]
	\arrow["{\mu_L}"', from=1-3, to=2-3]
\end{tikzcd}\]
By passing to the quotients by $G', P'$ and $L'$ in the respective columns we obtain
\[\begin{tikzcd}
	{X/G'=U/B'} & {X_P/P'=U/B'} & {X_L/L'=U_L/B_L'} \\
	{S/G'=GU/G'} & {S_P/P'=PU/P'} & {S_L/L'=LU_L/L'.}
	\arrow["{i'}"', from=1-2, to=1-1]
	\arrow["\mu", from=1-1, to=2-1]
	\arrow["{\mu_P}", from=1-2, to=2-2]
	\arrow["i"', from=2-2, to=2-1]
	\arrow["p", from=2-2, to=2-3]
	\arrow["{\mu_L}"', from=1-3, to=2-3]
	\arrow["{p'}", from=1-2, to=1-3]
\end{tikzcd}\]
The right square is Cartesian using $P/B\cong L/B_L.$ This implies that $\mu_*\un=i_*p^*\mu_{L,*}\un$ and we obtain the \emph{parabolic induction functor}
$$i_*p^*:\DKSpr(\mathfrak{S}_L)\to \DKSpr(\mathfrak{S}).$$
The corresponding maps between the Steinberg stacks
\[\begin{tikzcd}
	Z/G' & {Z_P}/P' & {Z_L}/L'
	\arrow["\pi", from=1-2, to=1-3]
	\arrow["\iota"', from=1-2, to=1-1]
\end{tikzcd}\]
yield an algebra homomorphism between the Hecke algebras
$$\iota_*\pi^*:\Hh(\mathfrak{S}_L)=(G_0(Z_L/L'),\star)\to \Hh(\mathfrak{S})=(G_0(Z/G'),\star).$$
Hence, we obtain the following theorem describing parabolic induction on the level of Springer $K$-motives.
\begin{theorem}\label{thm:parabolicinduction}
The following diagram commutes
\[\begin{tikzcd}
	{\DKSpr(\mathfrak{S}_L)} & {\DKSpr(\mathfrak{S})} \\
	{\Dperf(\Hh(\mathfrak{S}_L))} & {\Dperf(\Hh(\mathfrak{S})).}
	\arrow["{i_*p^*}", from=1-1, to=1-2]
	\arrow["{(\iota_*\pi^*)^*}", from=2-1, to=2-2]
	\arrow["\wr", from=1-1, to=2-1]
	\arrow["\wr", from=1-2, to=2-2]
\end{tikzcd}\]
\end{theorem}
\begin{proof} Recall that the vertical equivalences are constructed in \Cref{thm:propertyrgivesformality} from the weight complex functor for the weight structure generated by the objects $\mu_*\un$ and $\mu_{L,*}\un$, respectively. Since $\mu_*\un=i_*p^*\mu_{L,*}\un$, the parabolic induction functor $i_*p^*$ is weight exact. Now by \cite[Corollary 3.5]{sosniloTheoremHeartNegative2017} the weight complex functor commutes with weight exact functors. This implies the statement.
\end{proof}
\subsection{Specialization} Consider a Springer datum $\mathfrak{S}=(G,B,V,U,G\times \Gm).$ Abbreviate $\Hh=\Hh(\mathfrak{S})$ and $R=R(G\times \Gm)=K_0(\Bc (G\times\Gm))=K_0(\Bc G)\otimes \Z[\mathbf{q}^{\pm 1}]$ where $\mathbf{q}$ is a formal variable. We can identify
$$\Dperf(R)\stackrel{\sim}{\leftarrow}\DKconst(\Bc (G\times\Gm))=\langle\un \rangle_{stb, \inplus}\subset \DK(\Bc (G\times\Gm)).$$
Then $N\in \DKconst(\Bc (G\times\Gm))$ acts on $M\in \DKSpr(\mathfrak{S})$ via $M\otimes f^*(N)$ for $f: S/(G\times \Gm)\to \Bc (G\times \Gm).$ 

Similarly, $Z/(G\times\Gm)\to \Bc (G\times\Gm)$ yields a morphism of rings $R\to \Hh,$ whose image is central in $\Hh.$
This yields an action of $\Dperf(R)$ on $\Dperf(\Hh)$, which is compatible with the equivalence $\DKSpr(\mathfrak{S})\stackrel{\sim}{\to}\Dperf(\Hh).$

For an ideal $I\subset R$ we abbreviate $\Hh_I=\Hh\otimes_RR/I$ and denote by 
$$\DKSpr(\mathfrak{S})_I=\DKSpr(\mathfrak{S})\otimes_{\DKconst(\Bc (G\times\Gm))}\Dperf(R/I)$$
the Lurie tensor product. Then, we obtain the following diagram, describing modules over the specialized Hecke algebra in terms of Springer motives.
\[\begin{tikzcd}
	{\DKSpr(\mathfrak{S})} & {\DKSpr(\mathfrak{S})_I} \\
	{\Dperf(\Hh)} & {\Dperf(\Hh)_I} & {\Dperf(\Hh_I)}
	\arrow["\wr", from=1-1, to=2-1]
	\arrow["\wr", from=1-2, to=2-2]
	\arrow[from=2-1, to=2-2]
	\arrow[from=1-1, to=1-2]
	\arrow[from=1-2, to=2-3]
	\arrow["\sim"', from=2-3, to=2-2]
\end{tikzcd}\]
In particular, we will use this to specialize the formal variable $\mathbf{q}\in R$ to $q\in \Z,$ for example $q=p^n,$ and obtain the following.
\begin{theorem}\label{thm:specializeq} There is an equivalence
$$\DKSpr(\mathfrak{S})_{\mathbf{q}=q}\stackrel{\sim}{\to}\Dperf(\Hh_{\mathbf{q}=q}).$$
\end{theorem}
\subsection{Specialization and fixed points}
For maximal ideals, the specialization described in the previous section can also be obtained by passing to fixed points.  We explain this now, following \cite[Section 8.1]{chrissRepresentationTheoryComplex2010}. In this section, we assume that all objects and categories are scalar extended to the ground field $k$ which is algebraically closed of characteristic $0.$

Taking traces yields an identification
$$R\stackrel{\sim}{\to}\Oo((G^{ss}/\!/\op{Ad}(G))\times \mathbb{G}_{m})\stackrel{\sim}{\to} \Oo((T/\!/W)\times \mathbb{G}_{m})$$
where $G^{ss}$ denotes the semisimple elements in $G$ and $\op{Ad}(G)$ the adjoint action.

We obtain the natural map $R\to Z(\Hh)$  and $\Hh$ is a free $R$-module. For simplicity, we assume that the map $R\to Z(\Hh)$ is an isomorphism, which is for example the case for the affine Hecke algebra. 
Then, each point $a=(s,q)\in (G^{ss} \times \Gm)(k)$ yields a central character $\chi_{a}:Z(\Hh)\to L_a=k$ Denote by $\mathfrak{m}_a=\ker\chi_a$ and write $A\subset G\times \Gm$ for the closed subgroup generated by $a.$ 
Then, via \cite[(8.1.6)]{chrissRepresentationTheoryComplex2010} we obtain a chain of isomorphism of algebras
\begin{align*}
	\Hh_{\mathfrak{m}_a}&=\Hh\otimes_{Z(\Hh)}L_{a}\cong(G_0(Z/G'),\star)\otimes_{R}L_{a}\cong(G_0(Z/A),\star)\otimes_{K_0(BA)}L_{a}\\
	&\cong(G_0(Z^A),\star)\cong(\operatorname{CH}^\bullet(Z^A),\star).
\end{align*}
The last isomorphism is the Riemann--Roch map into the Chow-groups of $Z^A$. In particular, the algebra $\Hh_{\mathfrak{m}_a}$ admits a grading.

Without loss of generality, assume that $s\in T$ and denote by $G(s)\subset G$ and $W(s)\subset W$ the centralizers. Then, the fixed points $Z^A$ and $\mu^A: X^{A}\to S^{A}$ arise from another Springer datum $\mathfrak{S}_a$ for $G(s)$ indexed by $I=W/W(s).$ In particular, we can apply the formality results from \cref{thm:steinbergformality}.

We hence obtain the following categorical version of specialization.
\begin{theorem} There is a commutative diagram of equivalences
\[\begin{tikzcd}
	{\DKSpr(\mathfrak{S})_{\mathfrak{m}_a}} & {\DKSpr(\mathfrak{S}_a)} & {\DMSpr(\mathfrak{S}_a)} \\
	{\Dperf(\Hh_{\mathfrak{m}_a})} & {\Dperf(\Hh_{\mathfrak{m}_a})} & {\Dperf^\Z(\Hh_{\mathfrak{m}_a}).}
	\arrow["\wr", from=1-1, to=2-1]
	\arrow["v"', from=1-3, to=1-2]
	\arrow["\wr", from=1-3, to=2-3]
	\arrow["{v'}"', from=2-3, to=2-2]
	\arrow["\wr", from=1-2, to=2-2]
	\arrow[Rightarrow, no head, from=2-1, to=2-2]
	\arrow["\sim", from=1-1, to=1-2]
\end{tikzcd}\]
\end{theorem}
\begin{proof} The vertical equivalences are induced by the weight complex functors for the weight structures generated by $\mu_!\un$ and $\mu^A_!\un$, respectively. The top horizontal functors preserve these generators and are hence weight exact. By \cite[Corollary 3.5]{sosniloTheoremHeartNegative2017} the weight complex functor commutes with weight exact functors. This implies the statement. 
\end{proof}

\subsection{Classification of irreducibles}\label{sec:classificationofsimples}
Assume that $k=\overline{\Q}_\ell.$ Consider the $\ell$-adic realization functor $$\DMSpr_{\red}(\mathfrak{S}_a)\to \DSpr(\mathfrak{S}_a,\overline\Q_\ell).$$
The sheaf $\mu_!\overline\Q_\ell\in \DSpr(\mathfrak{S}_a,\overline{\Q}_\ell)$ is by the decomposition theorem \cite{beilinsonFaisceauxPervers1982} a direct sum of intersection cohomology complexes associated to irreducible local systems. These direct summands yield the irreducible representations of $\Hh_{\mathfrak{m}_a}$, see \cite[Theorem 8.6.12]{chrissRepresentationTheoryComplex2010}.
This yields the following classification of irreducibles of $\Hh$.
\begin{theorem}\label{thm:delignelanglands}
	There is a bijection between irreducible representations of $\Hh$ and pairs 
	$$\setbuild{(a,\mathcal{L})}{a\in (G^{ss}/\op{Ad}(G)\times \Gm)(k), \mathcal{L}\in \op{Irr}(\op{Loc}^{\op{Spr}}(S^A,\overline\Q_\ell))}$$
	where $\Ll$ is an isomorphism class of an irreducible local system that appears in the decomposition of $\mu_!\overline\Q_\ell\in \DSpr(\mathfrak{S}_a,\overline{\Q}_\ell).$
\end{theorem}
For the affine Hecke algebra, this yields the Deligne--Langlands parametrization, see \cite{kazhdanProofDeligneLanglandsConjecture1987}.

\section{Affine Hecke algebra and local Langlands}\label{sec:affineheckeandlocallanglands}
We discuss how the affine Hecke algebra and representations of split reductive $p$-adic groups arises from the setup discussed in \Cref{sec:kmotivicspringertheory}.
\subsection{Affine Hecke algebra}\label{sec:affinehecke} Let $G$ be a reductive group with simply-connected derived subgroup. Denote by $G\supset B\supset T$ a Borel subgroup and maximal torus. Let $W=N_G(T)/T$ the Weyl group, $X(T)=\Hom_{grp}(T,\Gm)=\Pic(BT)$ the character lattice and $W_{\op{aff}}=W\!\ltimes\!X(T)$ the affine Weyl group. Denote by $\Hh_{\op{aff},\mathbf{q}}$ the (generic) affine Hecke algebra whose underlying $\Z$-module is $\Z[\mathbf{q},\mathbf{q}^{-1}]\otimes\Z[W_{\op{aff}}]$ where $\mathbf{q}$ is a formal variable.

We consider the Springer datum, see \Cref{sec:springerdatum}, $$\mathfrak{S}=(G,B,\lgg^*,(\lgg/\lb)^*, G\times \Gm ).$$ Here, $\lgg\supset \lb$ are the Lie algebras of $G\supset B.$ Then $$\mu:\widetilde{\Nn}=G\times^B(\lgg/\lb)^*=T^*G/B\to\Nn\subset \lgg^*$$ is the Springer resolution of the nilpotent cone and $Z=\widetilde{\Nn}\times_\Nn\widetilde{\Nn}$ is the Steinberg variety. 
\begin{theorem}[\cite{kazhdanProofDeligneLanglandsConjecture1987,chrissRepresentationTheoryComplex2010, lusztigBasesEquivariantTheory1998}] There is an isomorphism of algebras $$\Hh(\mathfrak{S})=(G_0(Z/(G\times \Gm)),\star)\cong \Hh_{\op{aff},\mathbf{q}}.$$
\end{theorem}
We hence obtain from \Cref{thm:steinbergformality} the following $K$-motivic realization of the category of representations of the affine Hecke algebra.
\begin{theorem}\label{thm:formalityaffinehecke} There is an equivalence of categories
    $$\DKSpr(\Nn/(G\times \Gm))\cong\Dperf(\Hh_{\op{aff},\mathbf{q}}).$$
\end{theorem}
Moreover, we obtain the following integral version of the results of \cite{eberhardtSpringerMotives2021,eberhardtMotivicSpringerTheory2022} relating the representations of the \emph{graded} affine Hecke algebra to Springer motives.
\begin{theorem}\label{thm:formalitygradedaffinehecke} There is an equivalence of categories
    $$\DMSpr_{\red}(\Nn/(G\times \Gm))\cong\Dperf^\Z(\Hh^{\op{gr}}_{\op{aff},\mathbf{q}}).$$
\end{theorem}
We refer to \cite{riderFormalityNilpotentCone2013,antorFormalityDeligneLanglandsCorrespondence2023} for a similar result in the setting of $\ell$-adic sheaves.
\subsection{Representations of $p$-adic groups}
In the following, we assume that all objects and categories are $\C$-linear.
Denote by $\LD{G}$ the Langlands dual split reductive group over a non-Archimedean local field $F$ with finite residue field isomorphic to $\F_q.$ Let $I\subset \LD{G}$ be an Iwahori subgroup and $\Rep^I(G)$ the Abelian category of smooth $\LD{G}$-representations which are generated by their $I$-invariant vectors. Then, there is an isomorphism
$$C_c(I\backslash \LD{G}/I)\cong \Hh_{\op{aff},q}$$
between the convolution algebra of compactly supported $I$-biinvariant functions on $\LD{G}$ and the affine Hecke algebra $\Hh_{\op{aff},q}$ specialized at $\mathbf{q}=q.$ Passing to $I$-invariant vectors yields an equivalence
$$\Rep^I(\LD{G}(F))\stackrel{\sim}{\to}\Rep(\Hh_{\op{aff},q})$$
with the category of (finitely generated) representations of $\Hh_{\op{aff},q}.$ Using that the Hecke algebra has finite global dimension and \Cref{thm:specializeq} we obtain the following geometric realization of $\Rep^I(\LD{G}(F)).$
\begin{theorem}\label{thm:padicrepsviakmotives} There is an equivalence of categories
    $$\D^b(\Rep^I(\LD{G}(F)))\stackrel{\sim}{\to}\DKSpr(\Nn/(G\times \Gm))_{\mathbf{q}=q}.$$
\end{theorem}
Moreover, the results in \Cref{sec:parabolicinduction} and \Cref{sec:classificationofsimples} describe parabolic induction and the Deligne--Langlands parametrization of irreducibles in $\Rep^I(\LD{G})$ on the level of Springer $K$-motives.

\section{Quiver Hecke/Schur algebra} \label{sec:quiverheckeschuralgebra}
We discuss how ($K$-theoretic) quiver Hecke/Schur algebra arises from the setup discussed in \Cref{sec:kmotivicspringertheory}.
\subsection{Quiver representations}
Let $Q$ be a finite quiver with set of arrows, vertices, source and target function
$s,t: Q_1\rightrightarrows Q_0.$ Let $V=\bigoplus_{i\in Q_0}V_i$ a $Q_0$-graded vector space over $k$ with dimension vector $\dd=(\dd_i)=(\dim V_i)\in \Z_{\geq 0}^{Q_0}.$ We fix an identification $V_i=k^{\dd_i}.$ Denote by 
$$\Rep(V)=\Rep(\dd)=\bigoplus_{a\in Q_1}\Hom(V_{s(a)},V_{t(a)})$$
the space of representations of $Q$ with underlying vector space $V.$ Let $G_{\dd}=\GL(V)=\prod_{i\in Q}\GL(V_i).$ Then $G_\dd$ acts on $\Rep(\dd)$ by `conjugation'.

A \emph{composition} $\underline{\dd}=(\underline{\dd}^j)\in (\Z_{\geq 0}^{Q_0}-\{0\})^{\ell_{\dd}}$ of $\dd$ is a tuple of non-zero dimension vectors that sum to $\dd.$ A composition is called \emph{complete} if all parts $\underline{\dd}^j$ are unit vectors. We denote the sets of (complete) compositions of $\dd$ by $\op{Compf}(\dd)\subset \op{Comp}(\dd).$

 For $\underline{\dd}\in \op{Comp}(\dd),$ we denote by $$F_{\underline{\dd}}=(0=V^0\subset V^1\subset\dots\subset V^{\ell_d}=V)$$ 
the standard flag of $V$ such that $V^j/V^{j-1}$ has dimension vector $\underline{\dd}^j$ and by $P_{\underline{\dd}}\subset G$ the parabolic subgroup that stabilizes $F_{\underline{\dd}}.$ Moreover, denote by $\Rep(\dd)_{\underline{\dd}}\subset\Rep(\dd)$ the subspace of all representations $(\rho_a)_{a\in Q_1}$ that are compatible with the flag $F_{\underline{\dd}},$ that is, $\rho_a(V^j_{s(a)})\subset V^j_{t(a)}.$ Then $\Rep(\dd)_{\underline{\dd}}$ is stable under the action of $P_{\underline{\dd}}.$
\subsection{Hecke Algebra}
From this we obtain the \emph{quiver Hecke} and \emph{quiver Schur Springer datum}
\begin{align*}
    \mathfrak{S}_{\dd}^{\op{Hecke}}&=(\Rep(\dd),\op{Compf}(\dd),\{P_{\underline{\dd}}\},\{\Rep(\dd)_{\underline{\dd}}\},G_{\dd}\times \Gm)\text{ and}\\
    \mathfrak{S}_{\dd}^{\op{Schur}}&=(\Rep(\dd),\op{Comp}(\dd),\{P_{\underline{\dd}}\},\{\Rep(\dd)_{\underline{\dd}}\},G_{\dd}\times \Gm)
\end{align*}
indexed over all (complete) compositions, respectively.
The associated graded Hecke algebras
    $$\Hh^{\op{gr}}(\mathfrak{S}_{\dd}^{\op{Hecke}})=R_{\dd}\text{ and }
    \Hh^{\op{gr}}(\mathfrak{S}_{\dd}^{\op{Schur}})=A_{\dd}$$
are the \emph{quiver Hecke algebra} (also called \emph{KLR algebra}), see \cite{khovanovDiagrammaticApproachCategorification2009,rouquier2KacMoodyAlgebras2008,varagnoloCanonicalBasesKLRalgebras2011c}, and \emph{quiver Schur algebra}, see \cite{stroppelQuiverSchurAlgebras2014}, respectively.

In analogy, we write $$\Hh(\mathfrak{S}_{\dd}^{\op{Hecke}})=R^K_{\dd}\text{ and }
\Hh(\mathfrak{S}_{\dd}^{\op{Schur}})=A^K_{\dd}$$ for the ungraded versions, and call them the \emph{$K$-theoretic quiver Hecke algebra} and \emph{$K$-theoretic quiver Schur algebra}. It is an interesting question to describe these algebras diagrammatically, see \cite{zhouGeometryQuiverFlag2023}.
\begin{remark} Following the conventions for the affine Hecke algebra, see \Cref{sec:affinehecke}, one could argue that it would be more natural to call $\Hh^{\op{gr}}(\mathfrak{S}_{\dd}^{\op{Hecke}})=R_{\dd}$ the graded quiver Hecke algebra and $\Hh(\mathfrak{S}_{\dd}^{\op{Hecke}})=R^K_{\dd}$ simply the quiver Hecke algebra.
\end{remark}
So we obtain the following formality result for the $K$-theoretic quiver Hecke/Schur algebra.
\begin{theorem}
    There is an equivalence of categories
    \begin{align*}
        \DKSpr(\mathfrak{S}_{\dd}^{\op{Hecke}})&\stackrel{\sim}{\to}\Dperf(R^K_{\dd})\text{ and }\\
        \DKSpr(\mathfrak{S}_{\dd}^{\op{Schur}})&\stackrel{\sim}{\to}\Dperf(A^K_{\dd}).
    \end{align*}
\end{theorem}
Moreover, a similar result holds for the quiver Hecke/Schur algebra.
\begin{theorem}\label{thm:dmspringerforquiverheckeschur}
    There is an equivalence of categories
    \begin{align*}
        \DMSpr(\mathfrak{S}_{\dd}^{\op{Hecke}})&\stackrel{\sim}{\to}\Dperf(R_{\dd})\text{ and }\\
        \DMSpr(\mathfrak{S}_{\dd}^{\op{Schur}})&\stackrel{\sim}{\to}\Dperf(A_{\dd}).
    \end{align*}
\end{theorem}
In \cite{eberhardtMotivicSpringerTheory2022}, \Cref{thm:dmspringerforquiverheckeschur} was obtained in the special case of rational coefficients and quivers of type $ADE$ as well as the cycle quiver.

\printbibliography

\end{document}